\documentclass[12pt]{amsart}
\usepackage{latexsym,enumerate}

\usepackage{amsmath,amsthm,latexsym,amssymb,mathrsfs,xcolor}
\usepackage{comment}
\usepackage{graphicx}
\usepackage{cite}

\numberwithin{equation}{section}

\newtheorem{theorem}{Theorem}
\newtheorem{lemma}{Lemma}
\newtheorem{corollary}{Corollary}

\newtheorem{remark}{Remark}

\numberwithin{theorem}{section}
\numberwithin{corollary}{section}
\numberwithin{lemma}{section}
\numberwithin{definition}{section}
\numberwithin{proposition}{section}
\numberwithin{remark}{section}
\newcommand{\dint}{\displaystyle\int}

\newcommand{\R}{\mathbb{R}}

\title[]{Isoperimetric Bounds for Weighted Steklov Eigenvalues with Radial Weights}

\begin{document}

\author[F. Brock]{F. Brock$^1$}

\author[F. Chiacchio]{F. Chiacchio$^2$}

\setcounter{footnote}{1}
\footnotetext{
Martin-Luther-University of Halle, Landesstudienkolleg, 06114 Halle, 
Paracelsusstr. 22, Germany,  
e-mail: {\tt friedemann.brock@studienkolleg.uni-halle.de}
}

\vspace{.3 cm}

\setcounter{footnote}{2} \footnotetext{
Dipartimento di Matematica e Applicazioni ``R. Caccioppoli'', Universit\`a
degli Studi di Napoli Federico II, Complesso Monte S. Angelo, via Cintia,
80126 Napoli, Italy; e-mail: 
{\tt francesco.chiacchio@unina.it}
}

\begin{abstract} 
We study the following class of Steklov eigenvalue problems:
\[
\nabla \cdot \bigl( w \nabla u \bigr) = 0 \quad \text{in } \Omega, \qquad
\frac{\partial u}{\partial \nu} = \gamma v u \quad \text{on } \partial \Omega,
\]
where $w$ and $v$ are prescribed positive radial functions, $\Omega$ is a Lipschitz domain in $\mathbb{R}^N$ with $N \geq 2$ and $\nu$ denotes its outward unit normal. 
Extending classical results in the unweighted case due to Weinstock, the first author, and others, we establish isoperimetric inequalities for low-order eigenvalues under suitable symmetry assumptions on the domain.
In the first part, we consider the case $w(x) = |x|^{\alpha}$ and $v(x) = |x|^{\beta-\alpha}$, where the parameters $\alpha, \beta \in \mathbb{R}$ satisfy appropriate constraints. Our analysis relies on an explicit computation of the spectrum in the radial case, variational principles, and a family 
of weighted isoperimetric inequalities with ``double density''.
In the second part, we address the case $v \equiv 1$ and $w(x) = W(|x|)$, where $W$ is a non-decreasing, log-convex function. In this setting, the proof relies, among other tools, on a new weighted isoperimetric inequality, which may be of independent interest.
\\[0.2cm] 
\\[0.2cm]
\rm
\end{abstract}

\maketitle

\section{Introduction}
Consider the eigenvalue problem
\begin{equation}
\label{StekBc}
\left\{ 
\begin{array}{ll}
\Delta u = 0 & \text{in } \Omega,  \\[0.1cm] 
\displaystyle \frac{\partial u}{\partial \nu} = \gamma u & \text{on } \partial \Omega,
\end{array}
\right.   
\end{equation}
where $\Omega \subset \mathbb{R}^{N}$, $N\geq 2$, is a bounded Lipschitz domain and $\nu$ denotes its outward unit normal. This problem is usually referred to as the {\sl Steklov problem} (see \cite{Steklov}) and has attracted considerable attention over the last century; we refer the reader to \cite{Bandle}, \cite{CGGS}, \cite{FL}, \cite{Henrot}, \cite{H2} (Chapter 5), \cite{OV} and the references therein.

It is known that the ball maximizes the first non-trivial Steklov eigenvalue among all smooth sets $\Omega$ with fixed Lebesgue measure (see \cite{W} for the two-dimensional case and \cite{BrockSteklov} for arbitrary dimension). Moreover, the ball also minimizes the harmonic mean of the first $N$ non-trivial Steklov eigenvalues, as stated in the latter reference. In \cite{BDFRF}, the authors established a quantitative version of this inequality.

In this paper, we generalize these results to certain {\sl weighted} eigenvalue problems. More precisely, we consider
\begin{equation}
\label{weightedsteklov}
\left\{ 
\begin{array}{ll}
\nabla \cdot \bigl( w \nabla u \bigr) = 0 & \text{in } \Omega, \\[0.1cm]
\displaystyle \frac{\partial u}{\partial \nu} = \gamma v u & \text{on } \partial \Omega,
\end{array}
\right.
\end{equation} 
where $\Omega$ is a bounded Lipschitz domain in $\mathbb{R}^N$, $N\geq 2$, and $w$ and $v$ are given positive radial functions satisfying additional conditions to be specified later.

Standard theory for self-adjoint compact operators ensures that problem (\ref{weightedsteklov}) admits eigenfunctions $u_n$ with corresponding eigenvalues $\gamma_n = \gamma_n(\Omega)$ ($n \in \mathbb{N}_0 := \mathbb{N} \cup \{0\}$), with $\gamma_0 = 0 < \gamma_1 \leq \gamma_2 \leq \dots$, $u_0(x) \equiv 1$, and $\lim_{n \to \infty} \gamma_n = +\infty$.

Our approach is based on the use of suitable radial test functions in $\Omega$. Since the equations in (\ref{weightedsteklov}) are no longer translation-invariant, we assume that $\Omega$ satisfies one or both of the following symmetry properties:
\\[0.1cm]
${\bf (S_1)}$ 
  $\Omega$ 
is symmetric with respect to the origin, i.e.,
\begin{equation}
\label{sym1}
x \in \Omega \Longleftrightarrow -x \in \Omega.
\end{equation}
${\bf (S_2)}$ 
 $\Omega$ is invariant under rotations of $\pi/2$ about the origin in two-dimensional coordinate planes. That is, for any pair $(i,j)$ 
with $1\leq i<j \leq N$, we have
\begin{eqnarray}
\label{sym2}
&&
x= (x_1, \ldots , x_N ) \in \Omega \Longleftrightarrow
 y= (y_1 , \ldots , y_N) \in \Omega , \quad \mbox{where}  
\\
\nonumber 
&&
y_i = -x_j , \ y_j = x_i   \quad \mbox{and }\ y_k = x_k \ \mbox{for }\  k\in \{ 1, \ldots , N\} \setminus \{ i,j\} .
\end{eqnarray} 

\begin{figure}[h]
\centering
\includegraphics[scale=0.45]{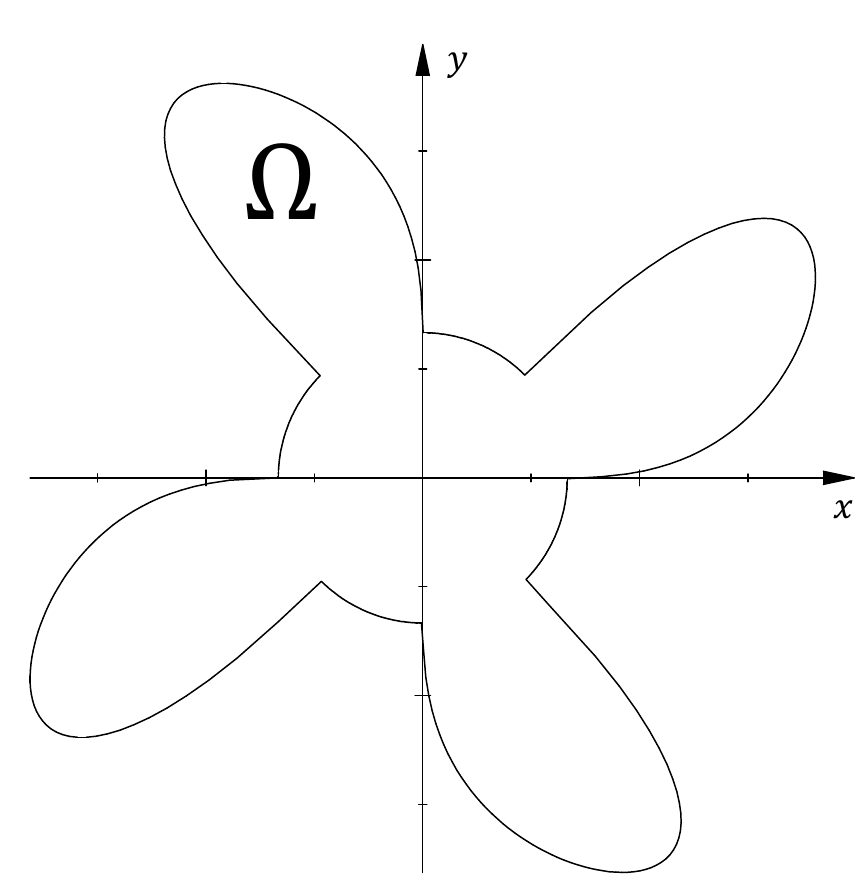}
\caption{A domain fulfilling condition ${\bf (S_2)}$}

\end{figure}

\begin{remark}
\leavevmode
\begin{enumerate}
    \item[\bf (a)]It is worth noting that condition \textbf{(S\(_2\))} implies \textbf{(S\(_1\))}  if $N$ is even.  Indeed, by applying a rotation of angle $\pi/2$ twice in any coordinate plane, $\Omega$ remains invariant under the simultaneous change of sign of the corresponding pair of coordinates. When $N$ is even, we may group all coordinates into disjoint pairs $(1,2), \dots, (N-1,N)$ to deduce that $\Omega$ is invariant under the map $x \mapsto -x$.
    \item[\bf (b)]Conversely, for odd $N$, \textbf{(S\(_2\))} does not imply \textbf{(S\(_1\))}. This is because the transformations generated by \textbf{(S\(_2\))} are orientation-preserving (having determinant $+1$), whereas the map $x \mapsto -x$ reverses orientation when $N$ is odd (having determinant $(-1)^N = -1$).
    \item[\bf (c)] The case where the origin lies on the boundary is excluded from this study, since under such a condition neither ${\bf (S_1)}$ nor ${\bf (S_2)}$ can be satisfied without violating the Lipschitz regularity of $\partial\Omega$. 
Indeed, assume by contradiction that $0 \in \partial\Omega$. 
 Since $\Omega$ is Lipschitz, after a rigid motion, we may assume that there exist $r>0$ and a Lipschitz function 
$$
\varphi:B_r'\to \mathbb R, \quad (B_r ^{\prime} : = \{ x' \in \mathbb{R} ^{N-1} :\, |x'|<r\} ),
$$
such that
\[
\varphi(0)=0
\ \mbox{ and }\ 
\Omega\cap Q_r
=
\{(x',x_N)\in Q_r:\ x_N>\varphi(x')\},
\]
where $Q_r:=B_r'\times (-r,r)$.
Then, for $t_0 >0$ sufficiently small, we have
\[
(0,\dots,0,t)\in \Omega \ \mbox{ for }\ 0<t<t_0.
\]
If ${\bf (S_1)}$ were satisfied, we would also have
\[
(0,\dots,0,-t)\in \Omega \ \mbox{ for } \ 0<t<t_0 ,
\]
which is impossible, since locally $\Omega$ lies only on one side of the graph $x_N=\varphi(x')$.

\noindent Likewise, if ${\bf (S_2)}$ were satisfied, then by applying twice a rotation of angle $\pi/2$ in a coordinate plane containing the $x_N$-axis, we would again obtain
\[
(0,\dots,0,-t)\in \Omega \ \mbox{ for }\ 0<t<t_0 ,
\]
which leads to the same contradiction.

\end{enumerate}
\end{remark}

Let us now state the main results of this paper. In the first part, we study the case where the weights are powers of the distance to the origin:
\begin{equation}
\label{alphabeta}
w = |x|^{\alpha}, \quad v = |x|^{\beta - \alpha}, \quad \alpha > -N, \ \beta \in \mathbb{R}.
\end{equation}

For convenience, we define
\begin{align}
\label{zdef}
z &:= \beta + 1 - \alpha, \\
\label{rhodef}
\rho &:= \sqrt{N^2 + 2\alpha (N-2) + \alpha^2}, \\
\label{elldef}
\ell &:= -N + \rho, \\
\label{kdef}
k &:= z + 1 - N + \rho,
\end{align}
and the function
\begin{equation}
\label{fdef}
f(z) := z(z+\rho)^2 + \rho \frac{(N-1)^2}{N}, \quad z \in \mathbb{R}.
\end{equation}

\medskip

\begin{theorem}\label{thm:1.1}
 Assume that $\Omega$ is a Lipschitz domain satisfying ${\bf (S_1)}$. Consider problem (\ref{weightedsteklov}) with $w$ and $v$ as in (\ref{alphabeta}). Suppose that one of the following conditions holds:
\begin{eqnarray*}
{\bf (i)} && z \geq 0;\\
{\bf (ii)} && -\frac{\rho}{N} \leq z \leq \min\{0, N-1-\rho\};\\
{\bf (iii)} && N \geq 3, \ N-1-\rho \leq z \leq 0, \mbox{ and if } 
N< \rho , 
\mbox{ then additionally } z \geq z_0,
\\
&&  \mbox{where $z_0$ is the only zero of $f$ in $[N-1-\rho,0]$};
\\
{\bf (iv)} && N=2, \quad \frac{|\alpha| - \rho}{2} \leq z \leq 0.
\end{eqnarray*}
Let $R>0$ be such that 
\begin{equation}
\label{StekBR}
\int_\Omega |x|^\ell \, dx = \int_{B_R} |x|^\ell \, dx.
\end{equation}
Then
\begin{equation}
\label{StekIsop}
\gamma_1(\Omega) \leq \gamma_1(B_R).
\end{equation}
\end{theorem}

Inequality (\ref{StekIsop}) can be improved under stronger symmetry assumptions.

\medskip
\begin{theorem}\label{thm:1.2}
\sl Assume the conditions of Theorem~\ref{thm:1.1} and additionally ${\bf (S_2)}$. Then
\begin{equation}
\sum_{i=1}^{N} \frac{1}{\gamma_i(\Omega)} \geq \frac{N}{\gamma_1(B_R)}.
\end{equation}
\end{theorem}

\medskip
\begin{corollary}
\label{cor:1.3}
  Theorems~\ref{thm:1.1} and \ref{thm:1.2} hold in the special case
\begin{equation}
\alpha = 0, \quad \beta \geq -2, \quad \ell = 0.
\end{equation}
\end{corollary}

In the Appendix, we provide  an example  showing that, without the symmetry assumption ${\bf (S_1)}$ for $\Omega$, the assertion of Theorems 1.1,
and a fortiori that of Theorem~\ref{thm:1.2},
may not hold. In other words, condition ${\bf (S_1)}$ is intrinsically related to the nature of the problem, rather than being a mere technical requirement for the proof to hold. 
See also Remark 4.3 in \cite{CdB}, where a similar phenomenon is addressed and clarified by means of a counterexample.

\medskip
The second part of the paper considers weights of the form
\begin{equation}
\label{Wlogconvex}
w(x) \in C^2(\mathbb{R}^N), \quad w(x) = W(|x|) > 0,
\end{equation}
with $W$ non-decreasing and log-convex (i.e., $\log W$ is convex), and
\begin{equation}
\label{v=1}
v \equiv 1.
\end{equation}

\medskip
\begin{theorem}\label{thm:1.4}
Assume that $\Omega$ is a Lipschitz domain satisfying ${\bf (S_1)}$.  Let $R>0$ be such that
\begin{equation}
\int_\Omega w \, dx = \int_{B_R} w \, dx.
\end{equation}
Then the inequality
\begin{equation}
\gamma_1(\Omega) \leq \gamma_1(B_R)
\end{equation}
holds.\end{theorem}

\medskip
\begin{theorem}
\label{thm:1.5}
 Assume the conditions of Theorem~\ref{thm:1.4} 
 and additionally ${\bf (S_2)}$.  Then
\begin{equation}
\sum_{i=1}^{N} \frac{1}{\gamma_i(\Omega)} \geq \frac{N}{\gamma_1(B_R)}
\end{equation}
holds.
\end{theorem}

\medskip
Now we outline the content of the article.

Section 2 deals with weights that are powers of the distance to the origin. First, we analyze (\ref{StekBc}) when $\Omega$ is a ball centered at the origin. In this case, the eigenfunctions are explicitly determined via separation of variables and then used as trial functions in variational characterizations of the first $N$ non-trivial eigenvalues of the problem in the original domain $\Omega$. Applying recent isoperimetric inequalities involving weighted perimeter and volume (Lemma~\ref{lem:2.1}) yields the eigenvalue estimates (Theorems~\ref{thm:1.1} and \ref{thm:1.2}).

Section 3 deals with weights satisfying (\ref{Wlogconvex}) and (\ref{v=1}). Again, radial trial functions are used in the variational characterization of eigenvalues in $\Omega$, and the desired estimates (Theorems 1.4 and 1.5) follow from a new weighted isoperimetric inequality (Lemma~\ref{lem:3.2}).
Finally, Section 4 (the Appendix) is devoted to the counterexample mentioned above.

\section{Steklov problems with weights $|x|^{\alpha }$}
\noindent
In this section we study the following weighted Steklov eigenvalue problem:
\begin{equation}
\label{Steklovxalpha} 
\left\{ 
\begin{array}{ll} 
\nabla \left( |x|^{\alpha } \nabla u \right) =0 & \mbox{ in } \ \Omega ,
\\[0.1cm]
\displaystyle{\frac{\partial u }{\partial \nu } = \gamma |x|^{\beta -\alpha } u } & \mbox{ on } \ \partial \Omega ,
\end{array}
\right.
\end{equation}
where
$\Omega$ is a bounded domain with Lipschitz boundary in $\R ^N $, $N\geq 2$,
$0\not\in \partial \Omega$, 
\begin{equation}
\label{alphabetaineq} 
\alpha >-N \ \mbox{ and }\  \beta \in \mathbb{R} .
\end{equation}

\noindent
{\bf 2.1. Preliminaries}
\\[0.1cm]
We denote by $L^p (\Omega , |x|^{\alpha })$ the weighted Lebesgue space of all measurable functions $u: \Omega \to \R $ with
$$
\left\Vert u\right\Vert _{L^{p}(\Omega ,\text{ }|x|^{\alpha })}
:= \left( \int_{\Omega }|u|^p   |x|^{\alpha } \, dx \right) ^{1/p} <\infty .
$$ 
The weighted Sobolev space $W^{1,p} (\Omega , |x|^{\alpha }  )$ is defined as the set of all functions $u\in L^p (\Omega , |x|^{\alpha }) $ having distributional derivatives $(\partial u / \partial x_i ) $, $i=1, \ldots ,N$, for which the norm
$$
\left\Vert u\right\Vert _{W^{1,p}(\Omega ,\text{ }|x|^{\alpha })}:=\left[
\left\Vert u\right\Vert _{L^{p}(\Omega ,\text{ }|x|^{\alpha
})}^{p}+\left\vert \left\Vert \nabla u\right\Vert \right\vert _{L^{p}(\Omega
,\text{ }|x|^{\alpha })}^{p}\right] ^{\frac{1}{p}}
$$
is finite.  
\\
\hspace*{0.3cm}
It is well-known (see, e.g., \cite{Necacs}) that, if $G$ is a bounded Lipschitz domain in $\mathbb{
R}^{N},$ with $N\geq 2$, then the trace operator $
W^{1,2}(G)\hookrightarrow L^{2}(G)$ is compact. Now the assumptions (\ref{alphabetaineq}) and $0\not\in \partial \Omega $ ensure that the embedding
of 
$W^{1,2}(\Omega ,|x|^{\alpha })$ 
into $L^{2}\left( \partial \Omega
,\,|x|^{\beta }\right) \equiv L^{2}\left( \partial \Omega \right) $ is
compact as well. 
This circumstance allows us to apply the classical spectral theory for
self-adjoint compact operators to problem (\ref{StekBc}).
\\
\hspace*{0.3cm}
Let  $\ell \in (-N,+\infty )$. 
We define the $\ell$-weighted measure of $\Omega$ as
\begin{equation}
\label{def_misura_l}
|\Omega|_{\ell} : = \int_{\Omega} |x|^{\ell} \, dx,
\end{equation}
and the $k$-weighted perimeter of $\Omega$ as
\begin{equation}
\label{def_perimetro_k}
P_{k}(\Omega) =
\begin{cases}
\displaystyle{\int_{\partial \Omega} |x|^k \, \mathcal{H}_{N-1} (dx) } & \text{ if } \Omega \text{ is } (N-1)\text{-rectifiable}\\[3ex]
+\infty & \text{otherwise}.
\end{cases}
\end{equation}
Here, $\mathcal{H}_{N-1}$ denotes the $(N-1)$-dimensional Hausdorff measure on $\mathbb{R}^N$.

In our analysis of the Steklov problem we will make use of the following
scale of isoperimetric inequalities in $\mathbb{R}^{N}$. The proof in the
cases \textbf{(i)'}-\textbf{(iii)'} below can be found  in \cite{ABCMP}, while
the case \textbf{(iv)'} was treated in \cite{McGillivray}. 
\\[0.1cm]
\begin{lemma}
\label{lem:2.1}
Let $N\in \mathbb{N}$, $\ell ,k\in \mathbb{R}$, $\ell >-N$, and
assume that one of the following conditions holds: \newline
{\bf (i)' } $N\geq 1$ and $k\geq \ell +1$; \newline
{\bf (ii)' } $N\geq 2$, $k\leq \ell+1$ and 
\begin{equation}
\label{k<0}
\ell \frac{N-1}{N}\leq k\leq 0;
\end{equation}
{\bf (iii)' } $N\geq 3$, $0\leq k\leq \ell+1$ and 
\begin{equation}
\label{caseiii}
\displaystyle{ \ell \leq \frac{(k+N-1)^{3}}{(k+N-1)^{2}-(N-1)^{2}/N}-N};
\end{equation}
{\bf (iv)' } $N=2$, $0\leq k\leq \ell +1$ and 
\begin{equation}
\label{N=2}
\ell\leq k-1+\frac{1}{k+1}.
\end{equation}
Then the following inequality holds for all smooth sets $\Omega \subset 
\mathbb{R}^{N}$: 
\begin{equation}
P_k (\Omega ) \geq C_{k,\ell,N}\left(
| \Omega | _{\ell} \right) ^{(k+N-1)/(\ell+N)},  
\label{isopRN}
\end{equation}
where 
\begin{equation}
C_{k,\ell ,N}:=(N\omega _{N})^{(\ell -k+1)/(\ell+N)}\cdot (\ell +N)^{(k+N-1)/(\ell +N)}.
\label{bestconst1}
\end{equation}
Equality in (\ref{isopRN}) holds for every ball in $\mathbb{R}^{N}$ centered
at the origin.
\end{lemma}

\begin{corollary} \label{cor:2.2}
Inequality~\eqref{isopRN} is equivalent to the following condition: if $|\Omega|_{\ell} = |B_R|_{\ell}$ for some $R>0$, then
\begin{equation} \label{isopRN2}
P_k (\Omega) \ge P_k (B_R).
\end{equation}
\end{corollary} 
For our application, we need to rewrite conditions \textbf{(i)'}--\textbf{(iv)'} in terms of the parameters $\alpha, \beta, z$, and $r$.

\begin{lemma}
\label{lem:2.3}
Let (\ref{zdef})--(\ref{kdef}) be in place and let the conditions {\bf (i)}--{\bf (iv)} be given as in Theorem~\ref{thm:1.1}. Then
\\
{\bf (i)} implies {\bf (i)'},
{\bf (ii)} implies {\bf (ii)'},
{\bf (iii)} implies {\bf (iii)'} and 
{\bf (iv)} implies {\bf (iv)'}. 
\end{lemma}

\noindent {\sl Proof: } 
 Recall that $\ell =-N+\rho$ and $k =z+1+\ell$.
\\
The first implication is obvious. 
\\
Assume that {\bf (ii)} holds. Then we have that $k\leq \ell +1$ and 
$$
(-N +\rho) \frac{N-1}{N} \leq z+1 -N +\rho \leq 0.
$$ 
The last two inequalities imply (\ref{k<0}).
Hence {\bf (ii)'} holds.
\\
Next assume {\bf (iii)}. From $N-1-\rho \leq z\leq 0$ we deduce $0\leq k\leq \ell +1$. 
\\
Now consider the function $f$ defined by (\ref{fdef}) for $z\in [N-1-\rho ,0]$. Note that 
\begin{equation}
\label{signoff}
f(N-1-\rho ) = \frac{(N-1)^3 }{N} \cdot (N-\rho ) \ \mbox{ and }\ 
f(0) = \rho \cdot \frac{(N-1)^2}{N} .
\end{equation}
$f' $ has 
has the two zeros $z_1 = -\rho$ and $z_2 =-\rho /3$ on $\mathbb{R}$, but at most the second one, $z_2 $, can be in $[N-1-\rho ,0]$.
 We split into two cases.
\\
{\sl (a)} Assume that $N\geq \rho $. Then $f(N-1-\rho )\geq 0$, and since $N\geq 3 $, we must have $N-1 -\rho \geq z_2 $. Hence $f'$ has no zeros on the interval $(N-1-\rho ,0)$. In view of (\ref{signoff}) we obtain
\begin{equation}
\label{f>0} 
f(z) \geq 0 \ \mbox{ on } \ [N-1-\rho ,0],
\end{equation}
and (\ref{caseiii}) follows.
\\
{\sl (b)} Assume that $N<\rho $, which implies $f(N-1-\rho) <0$. Since $f' $ has at most one zero in $[N-1-\rho ,0] $, $f$ has then exactly one zero, $z_0$, in $(N-1-\rho ,0)$ and there holds 
$$
f(z) >0 \ \mbox{ on }\ (z_0 ,0].
$$  
This implies (\ref{caseiii}) for these $z$.  
\\
We have shown {\bf (iii)'}.  
\\
Finally assume {\bf (iv)}.
Since $(|\alpha | -\rho)/2 \geq 1-\rho$, we have $1-\rho\leq z\leq 0$, which implies $0\leq k\leq \ell +1$.
Further, from $(|\alpha |-\rho)/2 \leq z$ we deduce $0\leq z + 1/(z+\rho)$, which implies $\ell \leq k-1 + 1/(k+1)$. This proves {\bf (iv)'}. 
$\hfill \Box $
\\[0.2cm]
{\bf 2.2. Formulation of the problem}
\\[0.1cm]
\noindent
For any function $u\in W^{1,2} (\Omega , |x|^{\alpha }) \setminus W_{0} ^{1,2} (\Omega , |x|^{\alpha })$ we define 
 the Rayleigh quotient 
\begin{equation}
R_{\Omega } (u) :=\frac 
{\dint_{\Omega }|\nabla u|^{2} \, |x|^{\alpha } \, dx}
{\dint_{\partial \Omega } u^{2} \, |x|^{\beta } \,  \mathcal{H}_{N-1} (dx)}.
\label{RayStek}
\end{equation}
The standard theory for
self-adjoint compact operators ensures that
 problem (\ref{Steklovxalpha}) has eigensolutions $u_{n}$ and the corresponding
eigenvalues $\gamma _{n}=\gamma _{n}(\Omega )$, 
($n \in \mathbb{N}_{0} := \mathbb{N}\cup \{0\}$), 
with $\gamma _{0}=0<\gamma _{1}\leq \gamma _{2}\leq \ldots $ , $
u_{0}(x)\equiv 1$, $\lim_{n\rightarrow \infty }\gamma _{n}=+\infty $ and 
\begin{equation}
\int_{\partial \Omega }u_{i}u_{j} \, |x|^{\beta } \mathcal{H}_{N-1} (dx) =0\text{ \
for }\ i\not=j,  \label{StekOrth1}
\end{equation}
and, equivalently, 
\begin{equation}
\int_{\Omega } \nabla u_{i} \cdot \nabla u_{j} \, |x|^{\alpha } \, dx=0\text{ \ for}
\ i\not=j.  \label{StekOrth2}
\end{equation}
Moreover we have the following variational characterizations (see \cite
{CourantHilbert}, \cite{Bandle} and \cite{HX}): 
\\[0.1cm]
{\bf 1.} For every $n\in \mathbb{N}$: 
\begin{eqnarray}
\label{StekVar1}
\gamma _{n}(\Omega )=R_{\Omega }(u_{n}) 
&=&
\displaystyle{ 
\inf 
\Bigg\{ 
R_{\Omega } (u) :\,  u\in W^{1,2} (\Omega , |x|^{\alpha } )\setminus W_0 ^{1,2} (\Omega , |x|^{\alpha } ), }
\\
\nonumber
&&  \qquad
\int_{\partial \Omega } u_{i}u \, |x|^{\beta } \, \mathcal{H}
_{N-1} (dx) =0,  \ \  i=0,\ldots ,n-1
\Bigg\} . 
\end{eqnarray}
{\bf 2.} For every $n\in \mathbb{N}$: 
\begin{eqnarray}
\label{StekVar2}
\sum_{i=1}^{n}\frac{1}{\gamma _{i}(\Omega )}
&=&
\sum_{i=1}^{n}\frac{1}{
R_{\Omega }(u_{i})}  
= 
\sup \Bigg\{ \sum_{i=1}^{n} \int_{\partial \Omega }v_{i}^{2}\,|x|^{\beta }\,
\mathcal{H}_{N-1} (dx) : 
\\
\nonumber
&& v_i \in W^{1,2} (\Omega , |x|^{\alpha }) \setminus W_0 ^{1,2} (\Omega , |x|^{\alpha })  , \, \int_{\partial \Omega }v_{i}\,|x|^{\beta }\, \mathcal{H}
_{N-1} (dx)=0,   \\
\nonumber 
&& 
\int_{\Omega }\,\nabla v_{i}\cdot \nabla
v_{j}\,|x|^{\alpha }\,dx=\delta _{ij},\ \ i,j=1,\ldots ,n \Bigg\}  . 
\end{eqnarray}
Let us first examine problem (\ref{Steklovxalpha}), when $\Omega $ is a ball centered at
the origin, i.e. $\Omega =B_{R}$, for some 
$R>0$, by means of separation of variables. We use $N$--dimensional
spherical coordinates $(r,\theta )$, where $r=|x|$ and $\theta =x|x|^{-1}\in 
\mathbb{S}^{N-1}$. 
The first equation in (\ref{Steklovxalpha}) can be rewritten as 
\begin{equation}
\left( r^{\alpha +N-1}u_{r}\right) _{r}+r^{\alpha +N-3}\Delta _{\theta }u=0,
\label{Sep1}
\end{equation}
where $\Delta _{\theta }$ denotes the Laplace-Beltrami operator on $\mathbb{S
}^{N-1}$. Setting $u(x)=Y(\theta )F(r)$, this becomes 
\begin{equation}
\frac{\left( r^{\alpha +N-1}F^{\prime }\right) ^{\prime }}{r^{\alpha +N-3}F}
=-\frac{\Delta _{\theta }Y}{Y}=:\overline{j}.  \label{Sep2}
\end{equation}
It is well known that (\ref{Sep2}) is fulfilled (see, e.g., \cite{Muller}) if and only if $
\overline{j}=j(j+N-2)$ with $j\in \mathbb{N}_0 = \mathbb{N} \cup \{0\}$. 
Hence (\ref{Sep2}) reduces to 
\begin{equation}
F^{\prime \prime }+(\alpha +N-1)\frac{F^{\prime }}{r}=j(j+N-2)\frac{F}{r^{2}}
\quad \mbox{on }\ (0,R),  \label{Sep3}
\end{equation}
which has the solutions $F(r)=r^{m}$, with 
\begin{equation}
m=m_{1,2}(j)=\frac{2-N-\alpha }{2}\pm \frac{1}{2}\sqrt{(N-2+\alpha
)^{2}+4j(j+N-2)},\quad (j\in \mathbb{N}_{0}).  
\label{defm}
\end{equation}
The case $j=0$ gives the eigenfunction $u_0 \equiv 1$ with corresponding eigenvalue $\gamma _0 = 0$, and  note that the orthogonality condition on $\partial B_R $ in (\ref{StekVar1}) ensures that this is also the only {\sl radial} eigenfunction.
\\ 
If $j\in \mathbb{N}$, then we have 
\begin{eqnarray*}
m_2 (j) &=& \frac{2-N-\alpha }{2} - \frac{1}{2}\sqrt{(N-2+\alpha
)^{2}+4j(j+N-2)} ,
\end{eqnarray*} 
and a short computation shows that the function $u(x) = Y(\theta ) r^{m_2 (j)} $ does not belong to the space 
$W^{1,2} (B_R , |x|^{\alpha})$.  
Hence we must have $m=m_{1}(j)(\geq
0) $, so that the solutions of (\ref{Steklovxalpha}) are of the form 
\begin{equation}
u=r^{m_{1}(j)}Y_{j}(\theta ),  \label{StekRad1}
\end{equation}
where $Y_{j}$ is an eigenfunction for the Laplace-Beltrami operator on $
\mathbb{S}^{N-1}$, with eigenvalue $\overline{j}=j(j+N-2)$, 
$(j\in \mathbb{N})$. \newline
Furthermore, from the second equation in (\ref{Steklovxalpha}) we
obtain 
$$
m_{1}(j)R^{m_{1}(j)-1}=\gamma R^{\beta -\alpha }R^{m_{1}(j)},
$$ 
therefore
$u$ is an eigenfunction in $B_{R}$ to the eigenvalue 
\begin{equation}
\gamma =m_{1}(j)R^{\alpha -\beta -1},\quad 
(j \in\mathbb{N}_{0}).
\label{StekRad2}
\end{equation}
From now on, let 
\begin{equation}
m:=m_{1}(1)=\frac{2-N-\alpha }{2}+\frac{1}{2}\sqrt{N^{2}+2\alpha
(N-2)+\alpha ^{2}}.  
\label{m=m1}
\end{equation}
Hence, the first  $N+1$ eigenfunctions and eigenvalues of problem
(\ref{Steklovxalpha}) for  $\Omega = B_{R}$, are given by
\begin{equation}
\label{StekRad0}
u_{0}(x) \equiv 1, \quad \gamma_{0}(B_{R}) = 0,
\end{equation}
and for $i=1,\ldots ,N,$
\begin{equation}
\label{StekRad4}
u_{i}(x) = x_{i} |x|^{m-1}, \quad \gamma_{i}(B_{R}) = m R^{\alpha - \beta - 1}.
\end{equation}
Let $i,j,k \in \{ 1, \ldots ,N\} $. From  (\ref{StekRad4})  we have 
\begin{equation*}
 \frac{\partial u_i }{\partial x_k} =
 \delta _{ik} |x|^{m-1} + (m-1) x_i x_k |x| ^{m-3} .
 \end{equation*}
 Hence 
\begin{eqnarray*}
\frac{\partial u_i }{\partial x_k} \cdot \frac{\partial u_j }{\partial x_k} &=& \delta _{ik} \delta _{jk} |x|^{2m-2} + (m-1) \delta _{jk} x_i x_k |x| ^{2m-4} + \\
&& + (m-1) \delta _{ik} x_j x_k |x|^{2m-4} + (m-1)^2 x_i x_j x_k ^2 |x|^{2m-6} .
\end{eqnarray*}
If $i\not= j$, this implies 
\begin{equation}
\label{productij}
\nabla u_i \cdot \nabla u_j = \sum_{k=1} ^N  \frac{\partial u_i }{\partial x_k} \cdot \frac{\partial u_j }{\partial x_k} 
= (m^2 -1) x_i x_j |x|^{2m-4} ,
\end{equation}
while for $i=j$ we obtain
\begin{equation}
\label{productii}
|\nabla u_i |^2 = \sum_{k=1} ^N \left(  \frac{\partial u_i }{\partial x_k} \right) ^2  = |x| ^{2m-2} + (m^2 -1) x_i ^2 |x|^{2m-4} .
\end{equation}
This finally yields
\begin{equation}
\label{sumiii}
\sum_{i=1} ^N |\nabla u_i |^2 = (N+ m^2 -1) |x|^{2m-2} .
\end{equation}
\noindent
{\bf 2.3. Proof of the isoperimetric results }
\\[0.1cm]
\noindent {\sl Proof of Theorem~\ref{thm:1.1}:} Let the functions $u_{i}$, ($i=0,1,\ldots ,N$), be given by
(\ref{StekRad0}) and (\ref{StekRad4}). In view of the symmetry assumption 
$\mathbf{(S_{1})}$ 
we have $\Omega = -\Omega := \{ -x:\, x\in \Omega \} $, which implies  
$$ \int_{\partial \Omega } x_i |x|^{m-1+\beta } \, \mathcal{H}_{N-1} (dx)= \int_{\partial \Omega } (-x_i ) |x| ^{m-1+\beta } \,   \mathcal{H}_{N-1} (dx) , \quad (i= 1, \ldots ,N).
$$
Hence
\begin{equation}
\int_{\partial \Omega }
 u_{i}(x) \, |x|^{\beta } \, \mathcal{H}_{N-1} (dx)
=\int_{\partial \Omega }
x_{i} |x|^{m-1+\beta } \, \mathcal{H}_{N-1} (dx)=0,\qquad
(i=1,\ldots ,N).  
\label{StekSym}
\end{equation}
Using the variational characterization (\ref{StekVar1}) for $\gamma
_{1}(\Omega )$, this means that 
\begin{equation}
\int_{\Omega }
|\nabla u_{i}|^{2} \,  |x|^{\alpha } \, dx
\geq 
\gamma _{1}(\Omega
)\int_{\partial \Omega }
 u_{i}^{2} \, |x|^{\beta } \, \mathcal{H}_{N-1} (dx).
\label{StekIneq1}
\end{equation}
In view of 
(\ref{productii}) this means 
\begin{equation}
\label{varinequality}
\int_{\Omega }\left( |x|^{\alpha +2m-2}+(m^{2}-1) \, x_{i}^{2} \, |x|^{\alpha
+2m-4}\right) \, dx\geq \gamma _{1}(\Omega )\int_{\partial \Omega
}x_{i}^{2} \, |x|^{\beta +2m-2} 
\, \mathcal{H}_{N-1} (dx),
\end{equation}
for $i=1,\ldots ,N$.
Defining $k$ and $\ell $ by 
(\ref{kdef}) and  (\ref{elldef}), we have
\begin{equation}
\label{kelldef2}
k = \beta +2m \ \mbox{ and }\ \ell = \alpha +2m -2.
\end{equation} 
Then adding (\ref{varinequality}) 
over all $i= 1, \ldots ,N$, and using 
(\ref{sumiii})
we obtain 
\begin{equation}
(N+m^{2}-1)\int_{\Omega }|x|^{\ell}\,dx\geq \gamma _{1}(\Omega )\int_{\partial
\Omega }|x|^{k}\,\mathcal{H}_{N-1} (dx).  \label{StekIneq3}
\end{equation}
Furthermore, repeating the above observations for $B_{R}$ in place of $
\Omega $, we find 
\begin{equation}
(N+m^{2}-1)\int_{B_{R}}|x|^{l}\,dx=\gamma _{1}(B_{R})\int_{\partial
B_{R}}|x|^{k}\, \mathcal{H}_{N-1} (dx).  \label{StekIneq4}
\end{equation}
Finally, Lemma~\ref{lem:2.1} and (\ref{StekBR}) yield 
\begin{equation}
 \label{StekIneq5}
\int_{\partial \Omega }|x|^{k}\, \mathcal{H}_{N-1}(dx) \geq \int_{\partial
B_{R}}|x|^{k}\, \mathcal{H}_{N-1} (dx), 
\end{equation}
Now 
the assertion (\ref{StekIsop}) follows from (\ref{StekBR}), (\ref{StekIneq3}), (\ref{StekIneq4}) and (\ref{StekIneq5}). 
$\hfill \Box $ 
\\[0.1cm]
\noindent {\sl Proof of Theorem~\ref{thm:1.2}:} 
Let $1\leq i<j\leq N$. Since $\Omega = \{ y= \varphi (x) :\, x\in \Omega \} $, where $x \longmapsto \varphi (x)=: y$ is the transformation used in condition $({\bf S_2)}$, we calculate
\begin{eqnarray*}
\int_{\Omega } x_i x_j |x|^{2m-4+\alpha } \, dx &=& \int_{\Omega } y_j (-y_i ) |y|^{2m-4+\alpha } \, dy  \quad \mbox{and }
\\
\int_{\Omega } \left( |x|^{2m-2+\alpha } + (m^2 -1) x_i ^2 |x| ^{2m-4+\alpha } \right) \, dx &=&  
\int_{\Omega } \left( |y|^{2m-2+\alpha } + (m^2 -1) y_j ^2 |x| ^{2m-4+\alpha } \right) \, dy.
\end{eqnarray*}
In view of  (\ref{productij}) and (\ref{productii}) this implies
\begin{eqnarray}
\label{S2a}
\int_{\Omega } \nabla u_i \cdot \nabla u_j |x|^{\alpha } \, dx &=& (m^2 -1) \int_{\Omega } x_i x_j |x|^{2m-4+\alpha } \, dx =0 \quad \mbox{and}
\\ 
\label{S2b}
\int_{\Omega }\left\vert \nabla u_{i}\right\vert ^{2}\left\vert
x\right\vert ^{\alpha }dx &=& \int_{\Omega }\left\vert \nabla u_{j}\right\vert ^{2}\left\vert
x\right\vert ^{\alpha }dx .
\end{eqnarray}
Since $i$ and $j$ were arbitrary, we obtain from this
\begin{eqnarray*}
\int_{\Omega }\left\vert \nabla u_{i}\right\vert ^{2} \left\vert
x\right\vert ^{\alpha }dxdy&=&\frac{N+m^{2}-1}{N}\int_{\Omega }\left\vert
x\right\vert ^{2m+\alpha -2}dx
\\&=&\frac{N+m^{2}-1}{N}\left\vert \Omega
\right\vert _{\ell}.
\end{eqnarray*}
Hence, if one defines 
\begin{equation*}
v_{i}:=\frac{u_{i}}{\sqrt{\dfrac{N+m^{2}-1}{N}\left\vert \Omega \right\vert
_{\ell}}}\text{ \ \ for }i=1,...,N,
\end{equation*}
(\ref{S2a}) and (\ref{S2b}) yield  
\begin{equation*}
\int_{\Omega }\nabla v_{i}\cdot \nabla v_{j}\text{ }\left\vert x\right\vert
^{\alpha }dxdy=\delta _{ij},\text{ \ }i,j=1,...,N.
\end{equation*}
Furthermore, the symmetry assumption ${\bf (S_1 )}$ gives 
$$
\int_{\partial \Omega } v_i |x|^{\beta } \, \mathcal{H}_{N-1} (dx) =  0, \quad i=1, \ldots ,N.
$$
Using the variational characterization (\ref{StekVar2}), we deduce that
\begin{eqnarray}
\label{Before_Isop}
\sum_{i=1}^{N}\frac{1}{\gamma _{i}(\Omega )}&\geq & \int_{\partial \Omega
}\left( \sum_{i=1}^{N}v_{i}^{2}\right) \left\vert x\right\vert
^{\beta }\text{ } \mathcal{H}_{N-1} (dx)
\\
\notag
&=&\frac{N}{\left( N+m^{2}-1\right)
\left\vert \Omega \right\vert _{\ell}}\int_{\partial \Omega }\left(
\sum_{i=1}^{N}u_{i}^{2}\right) \left\vert x\right\vert ^{\beta }
\text{ }\mathcal{H}_{N-1} (dx).
\end{eqnarray}
Since 
\begin{equation*}
\sum_{i=1}^{N}u_{i}^{2}=\sum_{i=1}^{N}x_{i}^{2}\left\vert x\right\vert
^{2m-2}=\left\vert x\right\vert ^{2m},
\end{equation*}
we can rewrite (\ref{Before_Isop})  as follows
\begin{eqnarray*}
\sum_{i=1}^{N}\frac{1}{\gamma _{i}(\Omega )}&\geq &
 \frac{N}{\left(
N+m^{2}-1\right) \left\vert \Omega \right\vert _{\ell}}
\int_{\partial \Omega }\left\vert x\right\vert ^{2m+\beta }\text{ } \mathcal{
H}_{N-1} (dx)
\\
 &=&\frac{N}{N+m^{2}-1} \cdot \frac{P_{k }\left( \Omega \right) }{
\left\vert \Omega \right\vert _{\ell }}.
\end{eqnarray*}
Finally, the weighted isoperimetric inequality, Corollary~\ref{cor:2.2}, gives the claim:
$$
\sum_{i=1}^{N}\frac{1}{\gamma _{i}(\Omega )}
\geq 
\frac{N}{N+m^{2}-1}
\frac{
P_{k }\left( B_{R}\right) }
{\left \vert B_{R} \right\vert _{\ell}}
=\sum_{i=1}^{N}\frac{1}{\gamma _{i}(B_{R})} = \frac{N}{\gamma _1 (B_R)}. \qquad \qquad \qquad \Box
$$
{\sl Proof of Corollary~\ref{cor:1.3}: } We make use of the conditions {\bf (i)'}, {\bf (iii)'} and {\bf (iv)'} of Lemma~\ref{lem:2.1}. Since $\alpha =0$, we have $m=1$, $\ell =0$ and $k= \beta +2 $. If $N\geq 2$ and $\beta \geq -1$, then {\bf (i)'} holds. Furthermore, if $N\geq 3 $ and $\beta \in [-2,-1]$, then {\bf (iii)'} is satisfied. Finally, if $N=2$ and $\beta \in [-2,1]$, then {\bf (iv)'} holds.  Now the assertion follows from Theorem~\ref{thm:1.1}, respectively Theorem~\ref{thm:1.2}. 
$\hfill \Box $ 


\section{Steklov problem with log-convex weight}
\noindent
In this section we study the following weighted Steklov problem:
\begin{equation}
\label{Steklovlogconvex} 
\left\{ 
\begin{array}{ll} 
\nabla \left( w \nabla u \right) =0 & \mbox{ in } \ \Omega ,
\\[0.1cm]
\displaystyle{\frac{\partial u }{\partial \nu } = \gamma  u } & \mbox{ on } \ \partial \Omega ,
\end{array}
\right.
\end{equation}
where
$\Omega$ is a bounded domain with Lipschitz boundary in $\R ^N $, $N\geq 2$, and $w$ is a positive radial weight satisfying (\ref{Wlogconvex}). 
\\[0.1cm]
{\bf 3.1. Preliminaries}
\\[0.1cm]
With the weight function $w$  we define a measure $\mu $ 
by
$$
d\mu := w(x)\, dx.
$$  
First we establish an Hardy-Littlewood-type inequality which is known to experts. For convenience we provide a simple proof.
\\[0.1cm]
\begin{lemma}\label{lem:3.1}
\sl (Hardy-Littlewood inequality) Let $H: [0,+\infty ) \to \mathbb{R}$ be measurable and  
$$
h (x):= H(|x|), \quad x\in \mathbb{R}^N .
$$
We assume that 
$$ 
r\longmapsto H(r) , \quad r\geq 0, 
$$
is non-decreasing. 
Furthermore, let $M \subset \mathbb{R}^N$ be measurable and $0<\mu (M )= \mu (B_R )<\infty $, ($R>0$). Then  
 \begin{equation}
 \label{HLinequality}
 \int_{B_R } h \, d\mu \leq \int_{M } h\, d\mu .
\end{equation}
\end{lemma}

\noindent {\sl Proof:} By our assumptions we have
\begin{eqnarray*}
&&
\mu (M \setminus B_R )  = \mu (B_R \setminus M ),
\\
&&
h (x)\geq H(R) \quad \mbox{in }\ M \setminus B_R \quad \mbox{and} 
\quad
h (x)\leq H(R) \quad \mbox{in }\ B_R \setminus M .
\end{eqnarray*}
This implies
\begin{eqnarray*}
\int_{B_R } h \, d\mu & = & \int _{B_R \cap M } h \, d\mu + \int _{B_R \setminus M } h \, d\mu
\\
& \leq &
\int _{B_R \cap M } h \, d\mu + H(R) \mu (B_R \setminus M )
\\
& \leq & \int _{B_R \cap M } h \, d\mu + \int _{M \setminus B_R  } h \, d\mu = \int_{M } h\, d\mu .
\qquad \qquad \qquad \Box 
\end{eqnarray*}
\hspace*{0.3cm}
Next we obtain a weighted isoperimetric inequality using the Divergence Theorem and the Hardy-Littlewood inequality. Notice that related results have been proven in \cite{Kolesnikov}, section 6.

\begin{lemma}
\label{lem:3.2}
 Let $W$ be given by (\ref{Wlogconvex}), and let $R>0$ be such that $\mu (\Omega )= \mu (B_R)$. Furthermore, let $\Phi \in C^2 [0, +\infty)$ be positive on $(0, +\infty )$ and such that
\begin{equation}
\label{monotoncond}
r\longmapsto \frac{N-1}{r} \Phi (r) + \Phi ' (r) + \frac{W'(r)}{W(r)} \Phi (r) \ \mbox{ is non-decreasing.}
\end{equation}
Then
\begin{equation}
\label{weightineq}
\int_{\partial B_R } w(x) \varphi (x) \, \mathcal{H}_{N-1} (dx) \leq 
\int_{\partial \Omega } w(x) \varphi (x) \, \mathcal{H}_{N-1} (dx),
\end{equation}
where $\varphi (x):= \Phi (|x|)$, ($x\in \mathbb{R} ^N $). 
\end{lemma}  

\noindent {\sl Proof:} From the Divergence Theorem we have 
($\nu$: exterior unit normal),
\begin{equation}
\label{divergencethm1}
\int_{\Omega } \mbox{div}\, \left( \frac{x}{|x|} w \varphi \right)\, dx = 
\int_{\partial \Omega } \frac{x\cdot \nu}{|x|} w\varphi \,   \mathcal{H}_{N-1} (dx) 
\leq 
\int_{\partial \Omega } w\varphi \,   \mathcal{H}_{N-1} (dx).
\end{equation}
In case of  $B_R $ in place of $\Omega $ we get the equality  
\begin{equation}
\label{divergencethm2}
\int_{B_R } \mbox{div}\, \left( \frac{x}{|x|} w \varphi \right)\, dx = 
\int_{\partial B_R } w\varphi \,   \mathcal{H}_{N-1} (dx).
\end{equation}
On the other hand, by assumption (\ref{monotoncond}) together with 
Lemma~\ref{lem:3.1} we obtain
\begin{eqnarray}
\label{HLineq}
\int_{\Omega } \mbox{div}\, 
\left( 
\frac{x}{|x|} w \varphi 
\right)\, dx
 &=& 
\int_{\Omega } 
\left( 
\frac{W'(|x|)}{W(|x|)} \Phi (|x|) +  \frac{N-1}{|x|} \Phi (|x|)+ \Phi ' (|x|) 
\right) \, d\mu 
\\
\nonumber
&\geq &
\int_{B_R } \left( 
\frac{W'(|x|)}{W(|x|)} \Phi (|x|) +  \frac{N-1}{|x|} \Phi (|x|)+ \Phi ' (|x|) 
\right) \, d\mu 
\\
\nonumber 
&= & 
\int_{B_R } \mbox{div}\, \left( 
\frac{x}{|x|} w \varphi 
\right) \, dx. 
\end{eqnarray}  
Now the assertion follows from (\ref{divergencethm1}), (\ref{divergencethm2}) and (\ref{HLineq}).
$\hfill \Box $

\begin{remark}
\label{rem:3.1}
Condition (\ref{monotoncond}) is satisfied for example if $\Phi (r) = r^m $ with $m\geq 1 $ and $W(r) = e^{V(r)}$, where $V$ is non-decreasing and convex.
\end{remark}

\noindent 
{\bf 3.2. Formulation of the problem}
\\[0.1cm]
Since $\Omega $ is bounded, we have
$$
C_1 \leq w(x)\leq  C_2 \quad \forall x \in \Omega 
$$
for two positive constants $C_1 $ and $C_2 $. Hence we may work in the classical (unweighted) functions spaces. We define the weighted Rayleigh quotient
$$ 
Q_{\Omega } (v) := 
\displaystyle{
\frac{\dint_{\Omega } |\nabla v|^2  \, d\mu }{\dint_{\partial \Omega } v^2 w\, \mathcal{H}_{N-1} (dx)}, \quad v\in W^{1,2} (\Omega ) \setminus W_0 ^{1,2} (\Omega ) .
}
$$
Analogously as in Section 2 we have:
Problem (\ref{Steklovlogconvex}) has eigensolutions $u_{n}$ and corresponding
eigenvalues $\gamma _{n}=\gamma _{n}(\Omega )$, 
($n \in \mathbb{N}_{0} $), 
with $\gamma _{0}=0<\gamma _{1}\leq \gamma _{2}\leq \ldots $ , $
u_{0}(x)\equiv 1$, $\lim_{n\rightarrow \infty }\gamma _{n}=+\infty $ and 
\begin{equation}
\label{orth1}
\int_{\partial \Omega }u_{i}u_{j} \, w \, \mathcal{H}_{N-1} (dx) =0 \ \mbox{
for }\ i\not=j,  
\end{equation}
or, equivalently, 
\begin{equation}
\label{orth2}
\int_{\Omega } \nabla u_{i} \cdot \nabla u_{j} \, d\mu =0 \ \mbox{  for }
\ i\not=j,  
\end{equation}
and there hold the following variational characterizations: 
\\[0.1cm]
{\bf 1.} For every $n\in \mathbb{N}$: 
\begin{eqnarray}
\label{Stekchar1}
\gamma _{n}(\Omega )=Q_{\Omega }(u_{n}) 
&=&
\displaystyle{ 
\inf 
\Bigg\{ 
Q_{\Omega } (u) :\,  u\in W^{1,2} (\Omega )\setminus W_0 ^{1,2} (\Omega ) , }
\\
\nonumber
&&  \qquad
\int_{\partial \Omega } u_{i}u \, w \, \mathcal{H}
_{N-1} (dx) =0,  \ \  i=0,\ldots ,n-1
\Bigg\} . 
\end{eqnarray}
{\bf 2.} For every $n\in \mathbb{N}$: 
\begin{eqnarray}
\label{Stekchar2}
\sum_{i=1}^{n}\frac{1}{\gamma _{i}(\Omega )}
&=&
\sum_{i=1}^{n}\frac{1}{
Q_{\Omega }(u_{i})}  
= 
\sup \Bigg\{ \sum_{i=1}^{n} \int_{\partial \Omega }v_{i}^{2}\, \,
\mathcal{H}_{N-1} (dx) : 
\\
\nonumber
&& v_i \in W^{1,2} (\Omega ) \setminus W_0 ^{1,2} (\Omega )  , \, \int_{\partial \Omega }v_{i}\, w \, \mathcal{H}
_{N-1} (dx)=0,   \\
\nonumber 
&& 
\int_{\Omega }\,\nabla v_{i}\cdot \nabla
v_{j} \, d\mu =\delta _{ij},\ \ i,j=1,\ldots ,n \Bigg\}  . 
\end{eqnarray} 
\hspace*{0.3cm}First we study the radial case, $\Omega =B_R $, by means of separation of variables. 
For convenience, we will often write $F$ and $W$, for $F(|x|)$ and $W(|x|)$, respectively.
\\[0.1cm]
Using $N$--dimensional
spherical coordinates $(r,\theta )$, 
the first equation in (\ref{Steklovlogconvex}) can be rewritten as 
\begin{equation}
\left( r^{N-1} Wu_{r}\right) _{r}+r^{N-3}W\Delta _{\theta }u=0.
\label{SepW1}
\end{equation}
Setting $u(x)=Y(\theta )F(r)$, this becomes 
\begin{equation}
\frac{\left( r^{N-1}W F^{\prime }\right) ^{\prime }}{r^{N-3} W F}
=-\frac{\Delta _{\theta }Y}{Y}=:\overline{j}.  \label{SepW2}
\end{equation}
As in Section 2 we have 
$\overline{j}=j(j+N-2)$ with $j\in \mathbb{N}_0 $. 
Hence (\ref{SepW2}) reduces to 
\begin{equation}
F^{\prime \prime }+\frac{W'}{W} F^{\prime} + (N-1) \frac{F^{\prime }}{r}=j(j+N-2)\frac{F}{r^{2}}
\quad \mbox{on }\ (0,R).  
\label{SepW3}
\end{equation}
The case $j=0$ gives the eigenfunction $u_0 \equiv 1$ with corresponding eigenvalue $\gamma _0 = 0$, and  note that the orthogonality condition on $\partial B_R $ in (\ref{Stekchar1}) ensures that this is also the only {\sl radial} eigenfunction. 
\\
For $j=1$ we obtain 
the first $N$ non-trivial eigenfunctions $u_1 , \ldots ,u_N $ of the problem. They have the form 
\begin{equation}
\label{ui}
u_i (x)= \frac{x_i}{|x|} F(|x|) , \quad x\in \mathbb{R} ^N , \ i= 1,\ldots ,N,
\end{equation}
where
\begin{eqnarray}
\label{F>0}
&&
\label{Fsmooth}
F\in C^2 [0, +\infty ), \quad
 F(0) = F'(0) =0 <F(r) \ \mbox{ for $r>0$, } 
\\
\label{odeW}
&& F^{\prime \prime} + \frac{W'}{W} F' + (N-1) \frac{F'}{r} = (N-1) \frac{F}{r^2 } \ \mbox{ on } (0,+\infty ), 
\end{eqnarray}
and there holds 
\begin{equation}
\label{gamma1N}
\gamma _1 (B_R ) = \ldots = \gamma _N (B_R ) >0 . 
\end{equation}
Note also that the function $F$ is unique up to constant multiples.
\\
\hspace*{0.3cm}Let $i,j,k \in \{ 1, \ldots ,N\} $. From  (\ref{ui})  we have 
\begin{equation}
\label{sumui}
\sum_{i=1} ^N (u_i )^2 = F^2 
\end{equation}
and 
\begin{equation*}
 \frac{\partial u_i }{\partial x_k} =
 \frac{\delta _{ik}}{|x|} F + \frac{x_i x_k }{|x|^2 } \left( F' - \frac{F}{|x|} \right) .
 \end{equation*}
 Hence 
\begin{eqnarray*}
\frac{\partial u_i }{\partial x_k} \cdot \frac{\partial u_j }{\partial x_k} &=& \delta _{ik} \delta _{jk} \frac{F^2}{|x|^2 } 
+ \frac{x_i x_j x_k ^2 }{|x|^4 } \left( F' - \frac{F}{|x|} \right)^2 +
\\
&& + \left( \delta _{ik} \frac{x_j x_k}{|x|^2} + \delta _{jk} \frac{x_i x_k }{|x|^2 } \right) \frac{F}{|x|} \left( F' - \frac{F}{|x|} \right).  
\end{eqnarray*}
If $i\not= j$, this implies 
\begin{equation}
\label{productijF}
\nabla u_i \cdot \nabla u_j = \sum_{k=1} ^N  \frac{\partial u_i }{\partial x_k} \cdot \frac{\partial u_j }{\partial x_k} 
= \frac{x_i x_j }{|x|^2 } \left( (F')^2 -\frac{F^2}{|x|^2 } \right),
\end{equation}
while for $i=j$ we obtain
\begin{equation}
\label{productiiF}
|\nabla u_i |^2 = \sum_{k=1} ^N \left(  \frac{\partial u_i }{\partial x_k} \right) ^2  = 
\frac{F^2}{|x|^2 } + \frac{x_i ^2 }{|x|^2} \left( (F')^2 - \frac{F^2}{|x|^2 } \right) .
\end{equation}
This also yields
\begin{equation}
\label{sumii}
\sum_{i=1} ^N |\nabla u_i |^2 = (F' )^2 + (N-1) \frac{F^2}{|x|^2} .
\end{equation}
\begin{lemma}\label{lem:3.3}
\sl Let $u_i$ be given by (\ref{ui})--(\ref{odeW}), $i=1,\ldots ,N$. Then
\begin{eqnarray}
\label{Fineq1}
&& \left( (F')^2 + 
\frac{N-1}{r^2} F^2 \right) ^{\prime} \leq   0  \quad \mbox{and}
\\
&&
\label{Fineq2}
\left( \frac{N-1}{ r} F^2 +2 F F' + \frac{W'}{W} F^2 \right) ^{\prime } \geq 0 \quad \forall r\in (0,R).
\end{eqnarray}
\end{lemma}
{\sl Proof:}  
Using (\ref{odeW}) we obtain 
\begin{eqnarray*}  
\left( (F')^2 + \frac{N-1}{ r^2 } F^2 \right) ^{\prime}
&=& 2F' F^{\prime \prime} -2 \frac{N-1}{r^3 } F^2 +2 \frac{N-1}{r^2 } FF'
\\
&=& -2 \frac{W'}{W} (F')^2 -2 \frac{N-1}{r } (F')^2 + 4\frac{N-1}{r^2 } FF' - 2 \frac{N-1}{r^3} F^2 
\\
&=& -2 \left( \frac{W'}{W} (F')^2 + \frac{N-1}{r} \left[F'- \frac{F}{r} \right] ^2 \right).   
\end{eqnarray*}
Since $W>0$ and $W' \geq 0$, (\ref{Fineq1}) follows.
\\
Finally, since $W$ is log-convex, we have that
$$
\left( \frac{W' (r)}{W(r)} \right) ' \geq 0 \quad \forall r\in [0, +\infty ) .
$$
Using this and once more (\ref{odeW}), we find 
\begin{eqnarray*}
&& \left( \frac{N-1}{r} F^2 +2 F F' + \frac{W'}{W}  F^2 \right) ^{\prime}
\\  
&=& 2 \frac{N-1}{r} FF' -\frac{N-1}{ r^2 } F^2 + 2(F^{\prime} )^2 + 2F F^{\prime \prime} + \left( \frac{W'}{W} \right) ' F^2  + 2 \frac{W'}{W} F F' 
\\
&=&
 \frac{N-1}{r^2} F^2 +2 (F')^2 +  \left( \frac{W'}{W} \right) ' F^2  \geq 0. \qquad \qquad \qquad \Box 
\end{eqnarray*}
{\bf 3.3. Proofs of the isoperimetric results}
\\[0.1cm]
{\sl Proof of Theorem~\ref{thm:1.4}: } Let $u_i $, $i=1,\ldots ,N$, be the first $N$ non-trivial eigenfunctions in $B_R $. 
\\
In view of assumption ${\bf (S_1 )}$ and (\ref{ui}) we have
\begin{equation}
\int_{\partial \Omega } u_i W \, \mathcal{H}_{N-1 } (dx) =0, \ \ i= 1, \ldots ,N.
\end{equation}
Hence the variational characterization of the eigenvalue $\gamma _1 (\Omega ) $ shows that 
$$
\int_{\Omega } |\nabla u_i | ^2 \, d\mu \geq \gamma _1 (\Omega ) \int_{\partial \Omega } u_i ^2 W \, \mathcal{H}_{N-1} (dx),   \quad i=1, \ldots ,N .
$$
Adding up and using the identities (\ref{sumui}) and (\ref{sumii}) we obtain
\begin{equation}
\label{EVineq1}
\int_{\Omega } \left( (F')^2 + \frac{N-1}{r^2 } F^2 \right) \, d\mu \geq \gamma _1 (\Omega ) \int_{\partial \Omega } F^2 W \, \mathcal{H}_{N-1} (dx) .
\end{equation}
Since we also have
$$
\int_{B_R } |\nabla u_i |^2 \,  d\mu = \gamma _1 (B_R ) \int_{\partial B_R} u_i ^2 W\, \mathcal{H}_{N-1} (dx), \quad i= 1, \ldots, N,
$$
we obtain 
\begin{equation}
\label{EVineq2}
\int_{B_R } \left( (F')^2 + \frac{N-1}{r^2 } F^2 \right) \, d\mu = \gamma _1 (B_R ) \int_{\partial B_R } F^2 W \, \mathcal{H}_{N-1} (dx) .
\end{equation}
Furthermore, Lemma~\ref{lem:3.1} and (\ref{Fineq1}) yield 
\begin{equation}
\label{HLappl} 
\int_{\Omega} \left( (F')^2 + \frac{N-1}{r^2 } F^2 \right) \, d\mu \leq
\int_{B_R} \left( (F')^2 + \frac{N-1}{r^2 } F^2 \right) \, d\mu .
\end{equation}
Putting $\Phi := F^2 $ and observing that (\ref{Fineq2}), we deduce (\ref{monotoncond}). Hence we may apply the isoperimetric inequality, Lemma~\ref{lem:3.2}, to obtain
\begin{equation}
\label{isopF2}
\int_{\partial B_R} F^2 W\, \mathcal{H}_{N-1} (dx) \leq 
\int_{\partial \Omega} F^2 W\, \mathcal{H}_{N-1} (dx).
\end{equation}
Now (\ref{StekIsop}) follows from 
(\ref{EVineq1})--(\ref{isopF2}). 
$\hfill \Box $.
\\[0.1cm]
{\sl Proof of Theorem~\ref{thm:1.5}:}  Using the assumption ${\bf (S_2 )}$, (\ref{productiiF}) and (\ref{sumii}) we find
\begin{equation}
\label{productiiF2}
\int_{\Omega } |\nabla u_i |^2 \, d\mu = \frac{1}{N} \int_{\Omega } 
\left(  (F' )^2 + (N-1 ) \frac{F^2 }{|x|^2 } \right) \, d\mu , \quad i= 1, \ldots ,N .
\end{equation} 
We normalize the function $F$ along with the eigenfunctions $u_i $ in $B_R$, such that
\begin{equation}
\label{normalize}
\int_{\Omega } |\nabla u_i |^2 \, d\mu =  1, \quad i= 1, \ldots ,N.
\end{equation}
Using the symmetry condition ${\bf (S_2) }$, (\ref{productijF}), (\ref{productiiF})
 and (\ref{normalize}), we obtain  
\begin{equation}
\label{scalarproduct}
\int_{\Omega }\nabla u_i \cdot \nabla u_j \, d\mu 
= \delta _{ij} \int_{\Omega } \frac{F^2 }{r^2 } \, d\mu + 
\int_{\Omega } \frac{x_i x_j }{r^2 }\left(  ( F')^2  - \frac{F^2}{r^2} \right) \, d\mu  = \delta _{ij} .\end{equation}  
Furthermore, there holds
\begin{equation}
\label{sumui2}
\sum_{i=1} ^N \int_{\partial \Omega} (u_i)^2 W \, \mathcal{H} _{N-1} (dx) = \int_{\partial \Omega } F^2 W\, \mathcal{H} _{N-1} (dx).
\end{equation} 
In view of (\ref{scalarproduct}) and (\ref{sumui2})
 we may use $ u_i $, $i=1, \ldots ,N$, as trial functions in (\ref{Stekchar2}) to obtain
\begin{equation}
\label{sumintegrals}
\sum_{i=1} ^N \frac{1}{\gamma _i  (\Omega )} 
 \geq  \sum_{i=1} ^N \int_{\partial \Omega } (u_i)^2 W\, \mathcal{H}_{N-1} (dx) = \int_{\partial \Omega } F^2 W \, \mathcal{H}_{N-1} (dx).
\end{equation}  
On the other hand,  
condition ${\bf (S_1 )}$ implies (\ref{isopF2}) (see the proof of Theorem~\ref{thm:1.4}). Hence
the assertion follows from (\ref{isopF2}), (\ref{sumintegrals}) and the identity
$$
\gamma_1 (B_R )\int_{\partial B_R } F^2 W \, \mathcal {H}_{N-1} (dx)= \int_{B_R } \left( (F')^2 + \frac{N-1}{r^2} F^2 \right) \, d\mu =N. \qquad \qquad \qquad \Box 
$$

\section{Appendix} 
Here we provide 
an example  showing that, without the symmetry assumption ${\bf (S_1)}$ for $\Omega$, the assertion of Theorem~\ref{thm:1.1},
and a fortiori that of Theorem~\ref{thm:1.2},
 may not hold.
\\[0.1cm]
\noindent 
 Consider the  case
\begin{equation*}
\alpha =0,\qquad \beta =-2,\qquad N\ge 2, \qquad l=0,
\end{equation*}
and, for $t>1$, let
\begin{equation*}
\Omega_t:=B_1(te_1)\subset \mathbb{R}^N,
\end{equation*}
where $e_1=(1,0,\dots,0)$. All these balls have the same volume $|B_1|$, where $B_1$ denotes the unit ball, $B_1(0)$, centered at the origin. Their first weighted Steklov eigenvalue is
\begin{equation*}
\gamma_1(\Omega_t)=\min\left\{\frac{\dint_{\Omega_t}|\nabla u|^2\,dx}{\dint_{\partial\Omega_t}\dfrac{u^2}{|x|^2}\,d\sigma_x}
:\ u\in H^1(\Omega_t)\setminus H^1_0(\Omega_t),\ \dint_{\partial\Omega_t}\frac{u}{|x|^2}\,d\sigma_x=0\right\}.
\end{equation*}
After the change of variables
\begin{equation*}
x=te_1+X,\qquad X\in B_1,
\end{equation*}
and setting
\begin{equation*}
v(X):=u(te_1+X),
\end{equation*}
this becomes
\begin{equation*}
\gamma_1(\Omega_t)=\min\left\{\frac{\dint_{B_1}|\nabla v(X)|^2\,dX}{\dint_{\partial B_1}\rho_t(X)\,v(X)^2\,d\sigma_X}
:\ v\in H^1(B_1)\setminus H^1_0(B_1),\ \dint_{\partial B_1}\rho_t(X)\,v(X)\,d\sigma_X=0\right\},
\end{equation*}
where
\begin{equation*}
\rho_t(X):=|te_1+X|^{-2},\qquad X\in \partial B_1.
\end{equation*}
Since $|X|=1$, one has
\begin{equation*}
t-1\le |te_1+X|\le t+1,
\end{equation*}
hence
\begin{equation*}
\frac{1}{(t+1)^2}\le \rho_t(X)\le \frac{1}{(t-1)^2}\qquad \forall X\in \partial B_1.
\end{equation*}
Therefore, for every admissible $v$,
\begin{equation*}
\int_{\partial B_1}\rho_t(X)\,v(X)^2\,d\sigma_X\le \frac{1}{(t-1)^2}\int_{\partial B_1}v(X)^2\,d\sigma_X,
\end{equation*}
and consequently
\begin{equation*}
\frac{\dint_{B_1}|\nabla v(X)|^2\,dX}{\dint_{\partial B_1}\rho_t(X)\,v(X)^2\,d\sigma_X}
\ge (t-1)^2\frac{\dint_{B_1}|\nabla v(X)|^2\,dX}{\dint_{\partial B_1}v(X)^2\,d\sigma_X}.
\end{equation*}
Hence
\begin{equation}
\gamma_1(\Omega_t)\ge (t-1)^2\inf\left\{\frac{\dint_{B_1}|\nabla v(X)|^2\,dX}{\dint_{\partial B_1}v(X)^2\,d\sigma_X}
:\ v\in H^1(B_1)\setminus H^1_0(B_1),\ \int_{\partial B_1}\rho_t(X)\,v(X)\,d\sigma_X=0\right\}.\label{eq:1}
\end{equation}
To estimate the infimum in \eqref{eq:1}, let
\begin{equation*}
m(v):=\frac{1}{|\partial B_1|}\int_{\partial B_1}v(X)\,d\sigma_X
\end{equation*}
be the Euclidean boundary mean of $v$. 
By construction, $v - m(v)$ has zero mean on $\partial B_1$; therefore, since the first nonzero Steklov eigenvalue of the unit ball is equal to $1$, it holds that
\begin{equation}
\int_{B_1}|\nabla v(X)|^2\,dX
=
\int_{B_1} |\nabla (v(X) - m(v))|^2 \, dX
\ge \int_{\partial B_1}(v(X)-m(v))^2\,d\sigma_X.
\label{eq:2}
\end{equation}
Now let
\begin{equation*}
\overline{\rho_t}:=\frac{1}{|\partial B_1|}\int_{\partial B_1}\rho_t(X)\,d\sigma_X.
\end{equation*}
Since $v$ is admissible, one has
\begin{equation*}
0=\int_{\partial B_1}\rho_t(X)\,v(X)\,d\sigma_X
=\overline{\rho_t}\int_{\partial B_1}v(X)\,d\sigma_X+\int_{\partial B_1}(\rho_t(X)-\overline{\rho_t})v(X)\,d\sigma_X.
\end{equation*}
Recalling that
\begin{equation*}
\int_{\partial B_1}v(X)\,d\sigma_X=|\partial B_1|\,m(v),
\end{equation*}
we obtain
\begin{equation*}
\overline{\rho_t}\,|\partial B_1|\,|m(v)|=\left|\int_{\partial B_1}(\rho_t(X)-\overline{\rho_t})v(X)\,d\sigma_X\right|.
\end{equation*}
The Cauchy--Schwarz inequality yields the estimate
\begin{equation*}
\overline{\rho_t} \, |\partial B_1| \, |m(v)| \le \|\rho_t - \overline{\rho_t}\|_{L^2(\partial B_1)} \, \|v\|_{L^2(\partial B_1)};
\end{equation*}
consequently, we obtain the following bound for $|m(v)|$
\begin{equation}
|m(v)| \le \frac{\|\rho_t - \overline{\rho_t}\|_{L^2(\partial B_1)}}{|\partial B_1| \, \overline{\rho_t}} \, \|v\|_{L^2(\partial B_1)}.
\label{eq:3}
\end{equation}
We now estimate the coefficient in \eqref{eq:3}. Since
\begin{equation*}
|te_1+X|^2=t^2+2tX_1+1,
\end{equation*}
we may write
\begin{equation*}
\rho_t(X)=\frac{1}{t^2}\frac{1}{1+\frac{2X_1}{t}+\frac{1}{t^2}}.
\end{equation*}
As $|X|=1$, one has
\begin{equation*}
\frac{2X_1}{t}+\frac{1}{t^2}=O(t^{-1})\qquad \text{uniformly on }\partial B_1.
\end{equation*}
Therefore
\begin{equation*}
\rho_t(X)=t^{-2}+O(t^{-3})\qquad \text{uniformly on }\partial B_1.
\end{equation*}
By averaging over $\partial B_1$, it follows that
\begin{equation*}
\overline{\rho_t} = t^{-2} + O(t^{-3}),
\end{equation*}
which implies
\begin{equation*}
\|\rho_t - \overline{\rho_t}\|_{L^{\infty}(\partial B_1)} = O(t^{-3}).
\end{equation*}
Consequently, we obtain
\begin{equation*}
\frac{\|\rho_t - \overline{\rho_t}\|_{L^2(\partial B_1)}}{\overline{\rho_t}} = O(t^{-1}).
\end{equation*}
Thus \eqref{eq:3} yields
\begin{equation*}
|\partial B_1|\,m(v)^2\le \frac{C}{t^2}\int_{\partial B_1}v(X)^2\,d\sigma_X,
\end{equation*}
and 
using the identity
\begin{equation*}
\int_{\partial B_1}(v(X)-m(v))^2\,d\sigma_X
=\int_{\partial B_1}v(X)^2\,d\sigma_X-|\partial B_1|\,m(v)^2,
\end{equation*}
we obtain
\begin{equation}
\int_{\partial B_1}(v(X)-m(v))^2\,d\sigma_X\ge \left(1-\frac{C}{t^2}\right)\int_{\partial B_1}v(X)^2\,d\sigma_X.
\label{eq:5}
\end{equation}
Combining \eqref{eq:2} and \eqref{eq:5}, we arrive at
\begin{equation*}
\int_{B_1}|\nabla v(X)|^2\,dX\ge \left(1-\frac{C}{t^2}\right)\int_{\partial B_1}v(X)^2\,d\sigma_X.
\end{equation*}
Substituting this into \eqref{eq:1}, we infer
\begin{equation}
\gamma_1(\Omega_t)\ge (t-1)^2\left(1-\frac{C}{t^2}\right).
\label{eq:6}
\end{equation}
In particular,
\begin{equation*}
\gamma_1(\Omega_t)\to +\infty\qquad \text{as }t\to\infty.
\end{equation*}
On the other hand, for the centered unit ball $B_1$ one has
\begin{equation*}
|x|^{-2}\equiv 1\qquad \text{on }\partial B_1,
\end{equation*}
hence
\begin{equation*}
\gamma_1(B_1)=1.
\end{equation*}
Therefore, for $t$ sufficiently large, we have
\begin{equation*}
\gamma_1(\Omega_t)>\gamma_1(B_1)
\end{equation*}
and consequently
\begin{equation*}
\sum_{i=1}^{N} \frac{1}{\gamma_i(\Omega_t)}
< 
\frac{N}{\gamma_1(B_1)}.
\end{equation*}

Hence the claim: the assertions of Theorems~\ref{thm:1.1} and \ref{thm:1.2} generally fail to hold if condition ${\bf (S_1)}$ is omitted.

\bigskip

\section*{Acknowledgements}
F. Chiacchio  was partially supported by Gruppo Nazionale per l'Analisi Matematica, la Probabilit\`a e le loro Applicazioni (GNAMPA) of Istituto Nazionale di Alta Matematica (INdAM) and by PNRR CUP E53D23018060001.

\end{document}